\magnification=1200
\baselineskip=18pt

\long\def\boxit#1{\leavevmode
\hbox{\vrule\vtop{\vbox{\kern.33333pt\hrule
    \kern1pt\hbox{\kern1pt\vbox{
\line{\hrulefill}\medskip 
#1}\kern1pt}}\kern1pt\hrule}\vrule}}

\font\eightpt=cmr8
\font\eightptit=cmti8
\font\eightptbf=cmbx8
\font\tenptit=cmti10
\font\tenptbf=cmbx10

\newcount\fnno\fnno=0

\long\def\fn#1{\advance\fnno by
1\def\it{\eightptit}\def\bf{\eightptbf}\ignorespaces
\footnote{$^{\hbox{\sevenrm\the\fnno}}$}
{\vtop{\advance\hsize by -20pt
                   \baselineskip=10pt\eightpt\noindent 
#1}}
\def\it{\tenptit}\def\bf{\tenptbf}     }

\font\small=cmr7
\overfullrule=0pt
\def\bull{\itemitem{$\bullet$}}

\def\Prob{\hbox{Prob}}
\def\os{\obeyspaces}
\def\C{\os\hbox{\bf C}\os}
\def\D{\os\hbox{\bf D}\os}

\def\s{\os\hbox{\bf s}\os}

\def\G{\os\hbox{\bf G}\os}

\def\L{\os\hbox{\bf L}\os}
\def\H{\os\hbox{\bf H}\os}
\def\HH{\os\hbox{\bf HH}\os}
\def\HT{\os\hbox{\bf HT}\os}
\def\TH{\os\hbox{\bf TH}\os}
\def\TT{\os\hbox{\bf TT}\os}
\def\CC{\os\hbox{\bf CC}\os}
\def\CD{\os\hbox{\bf CD}\os}
\def\DC{\os\hbox{\bf DC}\os}
\def\DD{\os\hbox{\bf DD}\os}
\def\T{\os\hbox{\bf T}\os}

\def\N{\os\hbox{\bf N}\os}
\def\P{\os\hbox{\bf P}\os}
\font\bigbold=cmbx12
\font\small=cmr8

\font\eightbf=cmbx8
\def\R{\os{\bf R}\os}
\def\Hom{\hbox{Hom}}
\def\Prob{\hbox{Prob}}

\centerline{\bigbold Quantum Strategies 
}
\centerline{by}
\centerline{Gordon B. Dahl}
\centerline{University of California, San Diego}
\centerline{and}
\centerline{Steven E. Landsburg}
\centerline{University of Rochester}\footnote{}{{\small 
We 
thank Mark Bils, Alan Stockman, Paulo Barelli, David 
Miller, 
David Meyer, Vince Crawford, Val Lambson, Navin 
Kartik and Joel Sobel for 
helpful 
remarks.}} 

\medskip

\medskip
\centerline{ABSTRACT}
\smallskip
{\narrower\narrower\eightpt\noindent
We investigate the consequences of allowing players to adopt
strategies which take advantage of quantum randomization devices. In games
of full information, the resulting equilibria are always correlated
equilibria, but not all correlated equilibria appear as quantum equilibria.
The classical and quantum theories diverge further in games of private
information. In the quantum context, we show that Kuhn's equivalence between
behavioral and mixed strategies breaks down. As a result, quantum technology
allows players to achieve outcomes that would not be achievable with any
classical technology short of direct communication; in particular they do
not occur as correlated equilibria.  
 In general, in games of
private information, quantum technology allows players to achieve outcomes
that are Pareto superior to any classical correlated equilibrium, but not
necessarily Pareto optimal. A simple economic example illustrates these
points.\par\par}

\bigskip

\noindent{\bf 1.  Introduction.}

Quantum game theory investigates the behavior of 
strategic agents with access to quantum technology.  
Such technology can be employed in randomization devices 
and/or in communication devices.  This paper is 
concerned strictly with quantum randomization.  For 
other recent papers on quantum randomization in economics
and game theory see, for example, [B] or [LaM].  For 
models of quantum communication see, for example, [EW], 
[EWL] or [La].

Just as a classical mixed strategy can be thought of as 
a strategy conditioned on the realization of some 
classical random variable, so a quantum strategy can be 
thought of as a strategy conditioned on the value of 
some quantum mechanical observable.  This is a non-
trivial expansion of the strategy space, because quantum 
observables are not bound by the laws of classical 
probability.  For example, if $X,Y,Z$ and $W$ are 
classical binary random variables, it is nearly trivial 
to prove that
$$\Prob(X\neq W)\le\Prob(X\neq Y)+\Prob(Y\neq 
Z)+\Prob(Z\neq 
W)\eqno(1.1)$$
The clear intuition is that the first and last elements 
of a sequence cannot differ unless two adjacent elements 
differ along the way.  Nevertheless, the existence of 
observables violating (1.1) is predicted by quantum 
mechanics and confirmed by laboratory experiments.  (The 
most glaring inconsistencies are avoided by the fact 
that neither $X$ and $Z$ nor $Y$ and $W$ can be observed 
simultaneously.) 

In Section 2 we will provide a simple informal example 
illustrating the failure of (1.1) and its consequences 
for game theory.  In the remainder of the paper, we will 
investigate the consequences of allowing players 
to adopt quantum strategies. As Levine ([Le]) has 
observed, the resulting equilibria (at least in games of 
full information) are always correlated equilibria in 
the sense of Aumann ([A]).  But not all correlated 
equilibria appear as quantum equilibria.  A correlated 
equilibrium $E$ might not be sustainable in a given 
quantum environment because no pair of quantum 
strategies yields the outcome $E$.  Moreover, even in 
the presence of such quantum strategies, $E$ might fail 
to be deviation-proof in the quantum context.   

The classical and quantum theories diverge further when 
we turn to games of private information.  Here, in the 
quantum context, Kuhn's equivalence between mixed and 
behavioral strategies ([K]) breaks down. 
 Of course a game of private information can always be 
modeled as a game of complete information with a more 
complex strategy space---but the point is that there is 
more than one way to construct that more complex 
strategy space.  Classically, the constructions are 
equivalent; not so in the quantum case.  As a result, 
quantum technology allows players to achieve outcomes 
that would not be achievable with any classical 
technology short of direct communication.

Players could in general do even better if they could 
condition their behavior on each others' signals as well 
as their own.   However, it's important to recognize 
that our quantum devices allow nothing of the kind.  
Neither player receives any information about the other 
player's signal.

As a general rule, in games of private information, 
quantum technology allows players to achieve outcomes 
that are Pareto superior to any classical correlated 
equilibrium, but  not necessarily Pareto optimal.
Examples in the final section will illustrate both 
points. 

A brief table of contents:

Section 2 offers a brief example to illustrate what we 
mean by the breakdown of (1.1); we will return to this 
simple example repeatedly throughout the paper.  

Section 3 establishes some notation and vocabulary 
for talking about classical equilibria (including mixed 
strategy and correlated equilibria) in a form that is 
suitable for generalization to the quantum context.
A game \G is imbedded in a larger game by giving players 
a choice of random variables on which to condition their 
strategies.  When each player is restricted to a single 
choice, we recover the notion of correlated equilibrium.  

Sections 4 and 5 provide the quantum generalization for 
games of complete information.  Here, instead of 
conditioning their strategies on classical random 
variables, players are permitted to condition their 
strategies on quantum observables.  We observe that (i) 
there are quantum strategies that no classical random 
variables can mimic and (ii) all quantum equilibria are 
also correlated equilibria (though not vice-versa).  
Section 4 explains the basic physics and Section 5 
applies it to game theory.

Sections 6 and 7 provide the classical and quantum 
generalizations for games where players receive private 
signals.  We demonstrate that Kuhn's mixed/behavioral 
strategy equivalence breaks down in quantum environments 
when there is private information.  In fact, while it is 
easy to define quantum games of behavioral strategies, 
we show that there is no reasonable way to define a 
quantum analogue of the Kuhnian games of mixed 
strategies.  We also show that with private information 
there can be quantum equilibria which are in no sense 
equivalent to classical correlated equilibria.  

In the remaining sections, we present some examples.  
Section 8 establishes some technical results that will 
be needed to compute the equilibria in those examples, 
and Section 9 establishes the examples themselves.  One 
technical result is deferred to the appendix.

\bigskip

\noindent {\bf 2.  Cats and Dogs.}

We begin with an informal example that illustrates the 
failure of (1.1) and its significance for game theory.  
The
analysis here is adapted from [CHTW].
We will return to this example in a more formal context 
in Section 7.

{\bf Example 2.1.}
Two players, who cannot communicate once the game is 
underway, are each asked one of two yes/no questions, 
e.g. ``Do you like dogs?'' or ``Do you like cats?''.  
Each player's question is chosen independently via a 
fair coin flip.  The players both win if and only if 
their answers agree, unless they both get the ``cats'' 
question, in which case they win if and only if their 
answers disagree.

In a classical environment, it's clear that players can 
achieve a Pareto optimal outcome by always agreeing, 
which yields a success rate of 3/4.  In particular, the 
players have nothing to gain by randomizing their 
responses.  

Now equip Player $i$ with a coin, which the player can 
observe after rotating it through either of two angles, 
$C$ or $D$.  The outcomes of 
these observations have the following probabilities:

$$\matrix
{\hbox{If both coins }&&
\hbox{If either coin}\cr
\hbox{are rotated through angle $C$}&&
\hbox{is rotated through angle $D$}\cr
\matrix{
&\hbox{\bf Coin Two}\cr
\matrix{\phantom{q}\cr
\matrix{ \vbox{\bf \hbox{Coin}\hbox{One}}&
\matrix{}\cr}\cr}&
\matrix{
&{\bf H}&{\bf T}\cr
\H&.15/2&.85/2\cr
\T&.85/2&.15/2 \cr}
\cr}
&\quad&
\matrix{
&\hbox{\bf Coin Two}\cr
\matrix{\phantom{q}\cr
\matrix{ \vbox{\bf \hbox{Coin}\hbox{One}}&
\matrix{}\cr}\cr}&
\matrix{
&{\bf H}&{\bf T}\cr
\H&.85/2&.15/2\cr
\T&.15/2&.85/2\cr}
\cr}\cr}\eqno(2.1.1)$$

\medskip

Given these coins, each player can adopt the following 
strategy:  

{\narrower\narrower\noindent\bf{Strategy 
2.1.2:}  \sl If I am asked the ``cat'' question, 
I will rotate my coin through angle $C$, and if I am 
asked 
the ``dog'' question I will rotate my coin through angle 
$D$.  
Either way, I will answer ``yes'' if and only if the 
coin shows heads.\par\par}  

A moment's reflection reveals that if both players adopt 
this strategy, they win 85\% of the time, which is a 
clear improvement over the classical maximum of 75\%.

Unfortunately for the players, no such coins exist in a 
world governed by the classical laws of physics and 
probability.  To see this, define the following random 
variables:

$X$ is the orientation (i.e. heads or tails) of Coin 
One, after rotation through Angle $C$.

$Y$ is the orientation of Coin Two, after rotation 
through Angle $D$.

$Z$ is the orientation of Coin One, after rotation 
through Angle $D$.

$W$ is the orientation of Coin Two, after rotation 
through Angle $C$.

Then chart (2.1.1) reveals that $\Prob(X\neq 
Y)=\Prob(Y\neq Z)=\Prob(Z\neq W)=.15$, while 
$\Prob(X\neq W)=.85$, so that (1.1) --- an easy theorem 
of classical probability theory --- is violated.  

Another way to say essentially the same thing is to note 
that no joint probability distribution for the random 
variables $X$, $Y$, $Z$ and $W$ can yield the values in 
Chart (2.1.1).

However, the laws of quantum mechanics do allow for the 
existence of such ``coins'' (actually, subatomic 
particles), which are routinely produced in physics 
laboratories and could plausibly be incorporated in the 
machinery of future quantum computers.

Therefore players equipped with quantum coins can do 
better than players equipped with classical random 
number generators.  On the other hand, they can't do 
arbitrarily well; in the example at hand, the laws of 
quantum mechanics set an upper limit of approximately 
85\%  (more precisely, $\cos^2(\pi/8)$) for the expected 
fraction of wins.

Just as the notion of a mixed strategy captures the 
options available to players equipped with 
independent random number generators, we will introduce 
the notion of a {\it quantum strategy\/} to capture the 
options available to players equipped with ``quantum 
coins''.  

{\bf Remark 2.2.}  
Readers unfamiliar with the relevant physics might be 
tempted to 
conclude that some sort of communication must take place 
between the players in Example 2.1.  Note however that 
no action by either player has any effect on the 
probability distribution of anything the other player 
can observe.  

In particular, beware the following fallacious argument:

{\narrower\narrower\noindent\bf{Argument 2.2.1.}\sl 
  Suppose both players have agreed to 
play Strategy 2.1.2.
Let $p$ denote the probability that Player One says 
``yes'' conditional on receiving the ``dog'' question.  
Then if Player Two receives the cat question, he knows 
that $p=.85$, and if Player Two receives the dog 
question, he knows that $p=.15$.  Therefore Player Two, 
on receiving his question, learns something 
previously known only to Player One, which implies some 
form of communication.\par\par}

One problem with this argument is that Player One never 
knows the value of $p$,  
Therefore nothing that Player 
Two learns about the value of $p$ can have been 
``communicated'' from Player One.  

By analogy, suppose that Players One and Two carry coins 
that are known to be the same color (because they were 
prepared that way by a referee).  When Player Two first 
observes his own coin, he immediately learns the color 
of Player One's coin --- but we daresay that nobody 
would want to call this an instance of communication.  
Careful reflection will convince the reader that the 
apparent ``communication'' in Example 2.1 is of exactly 
this nature.  Quantum technology overcomes classical 
restrictions on {\it correlations\/}, but all of the 
classical restrictions on {\it communication\/} remain 
intact.\fn{We note also that it is a well-established principle of 
physics that faster-than-light communication is 
impossible, whereas the apparent ``communication'' in 
this example is instantaneous.  It follows that the 
appearance of communication --- at least in any sense 
that a physicist would recognize --- must be deceptive.}

Finally, we remark that the effects of these quantum 
coins could certainly be mimicked by a mediator who 
observes both questions, observes Player One's strategy, 
and hands Player Two a weighted coin that has an 85\% 
chance of yielding the winning response.  But the whole 
point here is to model the capabilities and limits of 
quantum technology, not of mediators.  We can equally 
well imagine a mediator who coordinates the responses 
perfectly, but this is beyond the capability of our 
quantum coins.   

\bigskip

\noindent{\bf 3.  Classical Game Theory.}

In this section, we will review the basic notions of 
classical game theory (including correlated equilibria) 
in order to establish some (slightly nontraditional) 
notation and vocabulary that will be suitable for 
generalization to the quantum context.

Throughout this section we fix a two-player game \G with 
strategy sets $S_1,S_2$ and payoff functions $P_1,P_2$.  
For simplicity, we usually assume $S_1$ and $S_2$ are 
finite.

We also fix a probability space $\Omega$, which for 
concreteness we can take to be the unit interval.
A {\sl random variable\/} always means a random variable 
with domain $\Omega$.

In the traditional formulation of game theory, a mixed 
strategy is a probability distribution on the strategy 
space $S_i$.  It will be more convenient for us to think 
of a mixed strategy as a random variable with values in 
$S_i$.  Because a given probability distribution can be 
induced by many different random variables, this leads 
to a great proliferation of strategies, many of which 
are essentially interchangeable, and which we will want 
to think of as equivalent.  The defintions in 3.1 will 
clarify the notion of equivalence.

{\bf Definitions 3.1.}  Two stategies $s,t\in S_1$ are 
{\sl equivalent} if:
$$\hbox{For all $u\in S_2$ we have $P_1(s,u)=P(t,u)$ 
and $P_2(s,u)=P_2(t,u)$}\eqno(3.1.1)$$

Two strategies in $S_2$ are equivalent if they satisfy 
the obvious condition symmetric to (3.1.1).

We define the game $\overline{\G}$ by replacing $S_i$ 
with the set of all equivalence classes of strategies 
for player $i$ (and retaining the obvious payoff 
functions).  

We say that two games \G and \H are {\sl equivalent\/} 
if $\overline{\G}$ is isomorphic to $\overline{\H}$.

{\bf Definition 3.2.}  A {\it classical environment\/} 
is a pair $E=({\cal P}_1,{\cal P}_2$) where each ${\cal 
P}_i$ is a set of measurable partitions of $\Omega$.    

In what follows we fix a classical environment $E=({\cal 
P}_1,{\cal P}_2)$.  Let ${\cal X}_i$ be the set of 
$S_i$-valued random variables that are measurable with 
respect to some partition in ${\cal P}_i$.  

{\bf Definition 3.3.}  The game $\G(E)=\G({\cal 
P}_1,{\cal P}_2)$ is defined as follows:
\itemitem{\bull}Player $i$'s strategy set is ${\cal 
X}_i$.
\itemitem{\bull}Player $i$'s payoff function is 
$$P_i(X,Y)=\int_{S_1\times S_2}P_i(x,y)d\mu_{X,Y}(x,y)$$
where $\mu_{X,Y}$ is the probability distribution on 
$S_1\times S_2$ induced by $(X,Y)$. 

(We abuse notation slightly by using the same notation 
$P_i$ for the payoff functions in $\G$ and in $\G(E)$.)

We view the \G-strategy set $S_i$ as contained in the 
$\G(E)$-strategy set ${\cal X}_i$ by identifying $s\in 
S_i$ with the random variable that takes $\s$ as its 
only value.  We call these the {\it pure strategies\/}.

{\bf Example 3.4.}  Suppose that for each $i$, the 
random variables in ${\cal X}_i$ induce every possible 
probability distribution on the strategy set $S_i$.  
Suppose also  that every partition in ${\cal P}_1$ is 
independent from every partition in ${\cal P}_2$ (so 
that every random variable in ${\cal X}_1$ is 
independent of every random variable in ${\cal X}_2$).  
Then $\G(E)$ is equivalent (in the sense of 3.1) to the 
classical game of mixed strategies associated to $\G$.  

We will sometimes abuse language by calling $\G(E)$ {\it 
the\/} game of mixed strategies associated to \G, though 
the various choices for $E$ yield games that are only 
equivalent, not isomorphic.

{\bf Notation 3.5.} 
If ${\cal X}$ is any set of random variables, we write 
${\cal P}({\cal X})$ for the corresponding set of 
partitions of $\Omega$.   
If ${\cal X}_i$ is a set of $S_i$-valued random 
variables, we define the environment $E({\cal 
X}_1,{\cal X}_2)=({\cal P}_1,{\cal P}_2)$.  

{\bf Example 3.6.}  Let $X$ and $Y$ be $S_1$-valued and 
$S_2$-valued random variables.  It is immediate that the 
probability distribution induced on $S_1\times S_2$ by 
$(X,Y)$ is (in the usual sense) a correlated equilbirium 
in \G if and only if $(X,Y)$ is a Nash equilibrium in 
the game $\G(\{X\},\{Y\})$.  

{\bf Observation 3.7.}  Let $E$ be any classical 
environment for \G and suppose that $(X,Y)$ is a Nash 
equilibrium in the game $\G(E)$.  Then $(X,Y)$ is a 
correlated equilibrium in \G.

The converse to 3.7 does not hold:

{\bf Example 3.8.}  Let $S_1=S_2=\{\H,\T\}$.  

Let \G be the game with the following payoffs:\fn{In 
future examples, players will use coin flips to choose 
their strategies; therefore we've called the strategies 
\H and \T for ``heads'' and ``tails''.}

$$\matrix{
{\bf \hbox to 1.3in{}
\hbox{\bf Player Two}}\cr
{\vbox to 1in{\vfil \hbox{\bf Player
One}\vfil}
\hskip .2in\vbox {\offinterlineskip
\eightbf \halign{
\strut\hfil#&\quad\vrule#&\quad\hfil#
 \hfil&\quad\vrule#&
\quad\hfil#\hfil&\quad\strut\vrule#\cr
 &&{\bf H}&&{\bf T}&\cr
 \omit&&\omit&&\omit&\cr
 \noalign{\hrule}
\omit&&\omit&&\omit&\cr
   {\bf H}&&$(0,0)$&&$(2,1)$&\cr
 \omit&&\omit&&\omit&\cr
 \noalign{\hrule}
 \omit&&\omit&&\omit&\cr
 {\bf T}&&$(1,2)$&&$(0,0)$&\cr
 \omit&&\omit&&\omit&\cr
 \noalign{\hrule}
 }}}\cr}$$
\vskip 4pt\noindent

Let $X$, 
$Y$ and $W$ be binary random variables such that:

$$\Prob(X=W=\H)=\Prob(X=W=\T)=1/8\qquad
\Prob(X\neq W=\H)=\Prob(X\neq W=\T)=3/8$$
$$\Prob(Y=W=\H)=\Prob(Y=W=\T)=1/12\qquad \Prob(Y\neq
W=\H)=\Prob(Y\neq W=\T)=
5/12$$

Let $E=E(\{X,Y\},\{W\})$.

It is easy to check that both $(X,W)$ and $(Y,W)$ yield
correlated equilibria in \G.
But $(X,W)$ is not an equilibrium in the game $\G(E)$,
though $(Y,W)$ is.

{\bf Remark 3.9.}  Our modeling assumptions allow Player 
One to observe the realization of either $X$ or $Y$ but 
not both.  This is because ${\cal P}_1$ consists of two 
partitions.  An alternative model would define ${\cal 
P}_1$ to consist of a single partition that refines both 
of these, thus allowing Player One to play a strategy 
contingent on the realizations of both $X$ and $Y$.  
Depending on the intended real-world implementation, 
this might or might not be a better model.  

{\bf Remark 3.10.}  The moral of this section is this:
Every correlated equilibrium $(X,Y)$ occurs as the 
equilibrium in a game of the form $\G(E)$.  (In fact, we 
just take $E={\cal P}(X,Y)$ as defined in (3.4).)  
However, the same $(X,Y)$ might {\it not\/} be a 
correlated equilibrium in some other game $\G(E')$.  

This will become important when we get to the quantum 
analogue, because there (unlike here) physical 
considerations will dictate the choice of the quantum 
environment.  Thus, given a correlated equilibrium 
$(X,Y)$, the interesting question will not be ``Is 
$(X,Y)$ an equilibrium in {\it some\/} environment 
$E$?'' (to which the answer is trivially yes), but ``Is 
$(X,Y)$ an equilibrium in a {\it particular\/} 
environment $E$?'' (to which the answer is sometimes yes 
and sometimes no).

\bigskip

\noindent{\bf 4.  Quantum Measurements.}

Our next goal is to further expand players' options by 
allowing them to make quantum measurements.  In this 
section we will explain what that means, and in the next 
section we will incorporate these quantum measurements 
into our formal model.

{\bf Discussion 4.1:  States.}  
A classical coin is in one of two states, ``heads'' 
(which we denote \H) or ``tails'' (\T).  A ``quantum 
coin'' can be in any state of the form 
$$\alpha\H+\beta\T\eqno(4.1.1)$$
where $\alpha$ and $\beta$ are complex numbers, not both 
zero.  

We view (4.1.1) as an element of the two-dimensional 
complex vector space spanned by the symbols \H and \T.  
Two such vectors represent identical states if one is a 
(non-zero) scalar  multiple of the other.  

Thus, strictly speaking, a {\sl state\/} is not a vector 
but an equivalence class of vectors.  Nevertheless, we 
will often abuse language by using the same expression 
(4.1.1) for both a vector and the state that it 
represents.  

A heads/tails measurement of a coin in state (4.1.1) 
yields the outcome ``heads'' or ``tails'' with 
probabilities proportional to $|\alpha|^2$ and 
$|\beta|^2$.  When such a measurement is made, the state 
immediately changes to either \H or \T, depending on the 
measurement's outcome.

It goes without saying that no ordinary-sized coin obeys 
these laws of quantum physics, but spin-1/2 particles 
such as electrons do, with ``heads'' and ``tails'' 
replaced by ``spin up'' and ``spin down''.  For ease of 
comparison with the classical case, we will continue to 
speak of ``quantum coins''.

{\bf Discussion 4.2:  Transformations.}
Each physical action (such as rotating the coin through 
a pariticular angle) corresponds to some {\sl unitary\/} 
transformation $U$ of the state space.\fn{A unitary 
transformation is an invertible linear transformation 
$U$ such that $\overline{U}^T=U^{-1}$.  The overbar 
denotes complex conjugation.}  Performing the action 
transforms the penny's state from $\phi$ to 
$U\phi$.  

Because states are defined only up to multiplication by 
nonzero scalars, we can restrict attention to unitary 
transformations with determinant 1. These are called 
{\sl special unitary transformations\/} and are 
represented by matrices
$$\pmatrix{P&Q\cr -
\overline{Q}&\overline{P}\cr}\eqno(4.2.1)$$
where $P$ and $Q$ are complex numbers satisfying
$|P|^2+|Q|^2=1$.

Thus, for example, if a coin in state (4.1.1) is 
subjected to a physical operation corresponding to the 
transformation (4.2.1), then it is transformed into the 
new state
$$\alpha(P\H-\overline{Q}\T)+\beta(Q\H+\overline{P}\T)
=(\alpha P+\beta Q)\H+(-
\alpha\overline{Q}+\beta\overline{P})\T$$

{\bf Discussion 4.3:  Randomization}.
A quantum coin is neither more nor less useful than a 
classical randomizing device.  Given a coin in the state 
(4.1.1), a player can change the state at will by applying 
a unitary transformation.  Once the coin is in a state 
$\gamma\H+\delta\T$, it acts just like a weighted coin 
that comes up heads or tails with probabilities 
proportional to $|\gamma|^2$ and $|\delta|^2$.  All the 
new phenomena arise not from {\it single\/} coins but 
from interactions between {\it pairs\/} of coins.

{\bf Discussion 4.4:  Entanglement.}
Once a pair of coins have come into contact, they no 
longer occupy their own states.  Instead, the pair 
occupies a state
jointly represented by a non-zero vector
$$\alpha\HH+\beta\HT+\gamma\TH+\delta\TT\eqno(4.4.1)$$
As in the one-coin case, any non-zero multiple of this 
expression represents the same state.  Measurements 
yield the outcomes (heads,heads), (heads,tails), and so 
forth with probabilities proportional to $|\alpha|^2$, 
$|\beta|^2$ and so forth.  These probabilities hold even 
when the coins are examined at physically remote 
locations.

Coins that occupy a joint state of the form (4.4.1) are 
called {\sl entangled\/}.  

A physical operation on the first coin is represented by 
a unitary matrix which we can take to be of the form 
(4.2.1).  This transformation has the following effect 
on basis elements:
$$\HH\mapsto P\HH-\overline{Q}\TH$$
$$\HT\mapsto P\HT-\overline{Q}\TT$$
$$\TH\mapsto Q\HH+\overline{P}\TH$$
$$\TT\mapsto Q\HT+\overline{P}\TT$$
Its action on a general state of the form (4.4.1) is 
determined by these rules plus linearity.  Likewise, the 
same operation applied to the second coin has the 
following effect:
$$\HH\mapsto P\HH-\overline{Q}\HT$$
$$\HT\mapsto Q\HH+\overline{P}\HT$$
$$\TH\mapsto P\TH-\overline{Q}\TT$$
$$\TT\mapsto Q\TH+\overline{P}\TT$$

{\bf Notation 4.5.} Start with a coin in some state 
$\xi$.  Let Player One apply the transformation $U_1$ 
and let Player Two apply the transformation $U_2$.  The 
resulting state, computed according to the rules of 4.4, 
is denoted $(U\otimes 1)\xi(1\otimes V)$.  

(This notation will be familiar to readers familiar with 
the yoga of multilinear algebra; others can simply take 
it as a definition.)

{\bf Example 4.6.}  Suppose the two coins start in the 
{\sl maximally entangled state\/} \HH+\TT.  Players One 
and Two apply transformations $U$ and $V$ to the first 
and second coins.

The reader may check that the resulting state is given 
by
$$(U\otimes 1)\xi(1\otimes V)=\alpha\HH+\beta\HT+
\gamma\TH+\delta\TT$$
 where

$$UV^T=\pmatrix{\alpha&\beta\cr\gamma&\delta\cr}$$

In particular, suppose that
$$U=\pmatrix{\cos(\theta)&\sin(\theta)\cr
-sin(\theta)&\cos(\theta)\cr}\qquad\qquad
V=\pmatrix{\cos(\phi)&\sin(\phi)\cr
-sin(\phi)&\cos(\phi)\cr}$$
for some $\theta,\phi\in[0,2\pi)$.  (These 
transformations are implemented physically by rotating 
the coins through the angles $2\theta$ and $2\phi$.) 
Then the resulting state is
$$(U\otimes1)(\HH+\TT)(1\otimes V)=\cos(\theta-
\phi)\HH+\sin(\theta-\phi)\HT
-\sin(\theta-\phi)\TH+\cos(\theta-\phi)\TT$$
so that the resulting probability distribution over 
strategy pairs is
$$\Prob(\H,\H)=\Prob(\T,\T)=\cos^2(\theta-\phi)/2\qquad
  \Prob(\H,\T)=\Prob(\T,\H)=\sin^2(\theta-\phi)/2$$

Note that the unitary operators $U$ and $V$ affect the 
join probability distribution of the two coins, but not 
the individual probability distributions, so that each 
coin always turns up \H with probability 1/2.  More 
generally, the reader can check that regardless of the 
initial state, no choice of $U$ can affect the 
probability distribution of Player Two's outcomes (nor, 
of course, can the choice of $V$ affect the probability 
distribution for Player One).  This follows from the 
mathematics, or, if you prefer, from the physical 
principle that no influence can travel faster than 
light.  

Even if Player One chooses to ignore his coin (planning, 
perhaps, to play a pure strategy or a classical mixed 
strategy) then Player Two's coin reverts to a 
randomization device ala Discussion 4.3.

\bigskip

\noindent{\bf 5.  Quantum Environments.}

In this section we will formalize the options available 
to players 
with quantum coins.  We view the possession of a quantum 
coin as analogous to the ``possession'' of a set of 
measurable partitions as in Section 3.2.  That is, it 
leads naturally to an expansion of the player's strategy 
set.  

For ease of exposition, we will restrict ourselves to 
games \G in which the strategy spaces $S_1$ and $S_2$ 
each have cardinality two.  Everything generalizes 
easily to the case of arbitrarily large finite strategy 
spaces, and somewhat less easily to infinite strategy 
spaces.

{\bf Definition 5.1.}  Let \G be a two-player game with 
strategy spaces $S_1=S_2=\{\H,\T\}$.  Then a {\it 
quantum environment} for \G is a non-zero vector in the 
complex vector space spanned by $\HH,\HT,\TH,\TT$.  
($\xi$ represents the initial state of a pair of quantum 
coins, and is thus analogous to the pair of measurable 
partitions introduced in Section 3.)

In what follows, we fix a quantum environment 
$\xi$.

{\bf Notation and Conventions 5.2.}   We write ${\cal 
U}$ for the set of unitary transformations of the vector 
space spanned by \H and \T.  

{\bf Discussion 5.3.}  We want to define a game 
$\G(\xi)$ that expands each player's strategy space by 
allowing him to randomize over pure strategies via the 
selection of a unitary operator.  This will be the 
analogue of the game 
$\G(E)$ in Definition 3.3.  One small difference is that 
in the classical environment of 3.3, pure strategies are 
automatically included in the strategy sets ${\cal X}_i$
whereas they are not automatically included in the set 
${\cal U}$; therefore we have to append them.  Thus 
$\G(\xi)$ will be a game in which Player $i$'s strategy 
set is 
$${\cal U}_i={\cal U}\cup S_i\eqno(5.3.1)$$.  

The payoff functions are defined in the obvious way; a 
strategy pair imposes a probability distribution on 
$S_1\times S_2$ and we compute payoffs as expected 
values with respect to this probability distribution.  
The next two definitions will make this precise.

{\bf Definition 5.4.}  With notation as above, let 
$(U,V)\in {\cal U}_1\times {\cal U}_2$. 

We define a probability distribution $\mu_{UV}$ on 
$S_1\times S_2$ as follows:

\itemitem{a)}If $U$ and $V$ are both unitary operators,
write 
$$(U\otimes1)\xi(1\otimes V)=
\alpha\HH+\beta\HT+\gamma\TH+\delta\TT$$
(where the left hand side is defined in Section 4.5).  
Then the probabilities associated to $\HH, \HT$, and so 
forth are proportional to $|\alpha|^2, |\beta|^2$ and so 
forth.

\itemitem{b)}If $U$ is a unitary operator and $V=s$ is a 
pure strategy, write
$$(U\otimes1)\xi(1\otimes I)=
\alpha\HH+\beta\HT+\gamma\TH+\delta\TT$$
(where the $I$ on the left hand side is the identity 
transformation).
Then $\mu_{U,s}(\H,s)$ and $\mu(U,s)(\T,s)$ add to one 
and are proportional to the squared norms of the 
coefficients on $\H s$ and $\T s$.

\itemitem{c)}If $U=s$ is a pure strategy and $V$ is a 
unitary operator, we define $\mu_{UV}$ analogously to 
part b).

\itemitem{d)} If $U=s$ and $V=t$ are pure strategies, 
the probability distribution $\mu_{UV}$ is concentrated 
on $(s,t)$. 

{\bf Definition 5.5.}  Given a game \G and a 
quantum environment $\xi$, define a new game 
$\G(\xi)$ as follows:

\itemitem{a)}Player $i$'s strategy set is ${\cal 
U}_i={\cal U}\cup S_i$ (as in 5.3.1).  
\itemitem{b)}The payoff functions are given by 
$$P_i(X,Y)=\int_{S_1\times 
S_2}P_i(X_1,X_2)d\mu_{XY}(s,t)$$
where $\mu_(XY)$ is the probability distribution 
associated to $X$ and $Y$ according to Definition 5.4.

{\bf Example 5.6.}  Let \G be the game of Example 3.8 and 
let $\xi=\HH+\TT$.  

In the game $\G(\xi)$, suppose that
Player One chooses a unitary operator $U$.
Then, in light of Example 4.6, Player 2 clearly 
optimizes by choosing a unitary operator $V$ so that 
$$UV^T=\pmatrix{0&1\cr -1&0}\eqno(5.6.1)$$
(thereby putting all the probability weight on the 
outcomes \HT and \TH).  The same is true with the 
players reversed, so the equilibria are precisely the 
pairs $(U_1,U_2)$ satisfying (5.6.1).
In any equilibrium, the outcomes $(2,1)$ and $(1,2)$ are 
each realized with probability 1/2.  

It is easy to check that this is a correlated 
equilibrium, but it is more than that.  For example, the 
correlated equilibria $(X,W)$ and $(Y,W)$ of Example 3.8 
are not sustainable as quantum equilibria in this 
environment.

{\bf  Remarks 5.7.}  We want to single out those 
quantum environments that 
can be mimicked by classical technology.  

To that end, continue to let $S_i$ be Player $i$'s 
initial 
strategy set and let $\Omega$ be the unit interval, 
thought of as a probability sample space.  Let ${\cal 
F}_i$ be the set of all $S_i$-valued random variables 
defined on $\Omega$.  Given $(F,G)\in {\cal F}_1\times 
{\cal F}_2$, let $\mu_{FG}$ be the probability 
distribution on $S_1\times S_2$ induced by the random 
variable $F\times G$. 

We'll say that the quantum enviroment 
$\xi$ is {\sl classical}
if each strategy in ${\cal U}_i$ can be mimicked by some 
random variable in ${\cal F}_i$.  More precisely:

{\bf Definition 5.8.}  With notation as in 5.7, 
we say that $\xi$ is {\sl classical\/} if there are maps
$$\phi_i:{\cal U}_i\rightarrow{\cal F}_i$$ such that for 
any $(U,V)\in {\cal U}_1\times{\cal U}_2$, the  
probability distributions $\mu_{UV}$ 
and $\mu_{\phi_1(U)\phi_2(V)}$ are identical.

(Reminder: $\mu_{UV}$ is defined in 5.4.)  

{\bf Example 5.9.}  The quantum environment 
$\xi=\HH+\TT$ is not classical.  

To see this, define, for each $\theta\in[0,2\pi)$, the 
following unitary operator:
$$M(\theta)=\pmatrix{\cos(\theta)&\sin(\theta)\cr
-\sin(\theta)&\cos(\theta)\cr}$$

(Physically, the operator $M(\theta)$ corresponds to 
rotating the coin through the angle $2\theta$.)  

Construct the $\{\H,\T\}$-valued ``random variable''
$\hat{M}_i(\theta)$  by applying
the operation $M(\theta)$ to coin $i$ and then observing
the coin's orientation.

Now let
$$X={\hat{M}_1(-\pi/8)}\quad
          Y={\hat{M}_2(0)}\quad           
Z={\hat{M}_1(\pi/8)}\quad
            W={\hat{M}_2(\pi/4)}$$

Then by the computations in 4.6, we have
$$\matrix{\Prob(X\neq Y)=\Prob(Y \neq Z)=
\Prob(Z \neq W)= \sin^2(\pi/8)\approx .15\cr
\phantom{a}\cr
          \Prob(X\neq W)=\sin^2(3\pi/8)\approx .85\cr}$$
so that condition (1.1) is violated.  Thus no classical 
random variables $X,Y,Z,W$ can mimic the effects of the 
strategies $M(0),M(\pi/8),M(\pi/4),M(3\pi/8)$.  

{\bf Remark 5.10.}  It is not hard to prove (and will be 
immediately obvious to readers with appropriate physics 
backgrounds) that $\xi$ is classical if and only  if it 
is of the form 
$$\xi=\alpha\HH+\beta\HT+\gamma\TH+\delta\TT$$
where $\alpha\delta-\beta\gamma=0$.

{\bf Remark and Notation 5.11.}  A quantum environment 
$\xi$ expands players' strategy sets by allowing them to 
make quantum observations.  We can expand the strategy 
sets still further by allowing players to randomize over 
those observations:  If $E$ is a classical environment 
(as defined in 3.2) we can consider the game 
$\G(\xi)(E)$.  We will abbreviate this game by 
$\G(\xi,E)$.  

{\bf Example 5.12.}  Let \G be a two by two game and let 
$\xi$ be the quantum environment $\xi=\HH+\TT$.  Then if 
both players adopt quantum strategies, either player can 
force a probability distribution that is concentrated 
either on the main diagonal or the off-diagonal (with 
the the two boxes on that diagonal represented 
equiprobably).  (If Player One plays $U$, player Two can 
play $\overline{U}$ to force the main diagonal, or 
$\overline{U}\cdot\pmatrix{0&1\cr -1&0\cr}$ to force the 
off-diagonal.)  Unless the players agree on which 
diagonal is more desirable, this means there can be no 
equilibrium in pure quantum strategies, but there can 
certainly be equilibria in mixed quantum strategies, 
and/or correlated equilbria in quantum strategies i.e. 
equilibria in games of the form $\G(\xi,E)$.

\bigskip

\noindent{\bf 6. Games of Private Information.}

We begin this section with a slight reformulation of the 
classical theory of games with private information.  Our 
setup is clearly equivalent to the standard formulation, 
but better suited to quantum generalization.  

{\bf Definition 6.1.}  A (two-player) {\sl game of 
private information\/} consists of two {\sl strategy 
sets $S_i$\/}, two {\sl signal sets} ${\cal A}_i$, a 
probability distribution $\nu$ on ${\cal A}_1\times{\cal 
A}_2$, and two payoff functions
$$P_i:{\cal A}_1\times {\cal A}_2\times S_1\times 
S_2\longrightarrow \R$$
Unless stated otherwise, we will assume the strategy and 
signal sets are finite, though all this could be 
generalized.

We call these ``games of private information'' because 
all of our equilibrium concepts will assume (in effect) 
that Nature selects a signal pair $(a,b)\in {\cal 
A}_1\times{\cal A}_2$ and that each player knows only 
his own signal. 

{\bf Notation 6.2.}  For any sets $X,Y$, we write $\Hom(X,Y)$ 
for the set of all functions from $X$ to $Y$.

In the discussion to follow, we fix a game of private 
information $\G=(S_1,S_2,{\cal A}_1,{\cal 
A}_2,\nu,P_1,P_2)$.  

{\bf Definition 6.3.} The {\it associated game} $\G^\#$ 
has strategy sets $S_i^\#=\Hom({\cal A}_i,S_i)$ and 
payoff functions $$P_i^\#(F_1,F_2)=\int_{{\cal 
A}_1\times {\cal A}_2} 
P_i(A_1,A_2,F_1(A_1),F_2(A_2))d\nu$$

We think of $\G^\#$ as the game that results from 
\G when we introduce contingent strategies. 

A {\it Nash equilibrium\/} in \G is (by definition) a 
Nash equilibrium in $\G^\#$.

{\bf Reminders 6.4.}  Now we will allow players to 
randomize.  As in Section 3, we work with a fixed sample 
space $\Omega$.  Recall from Definition 3.2 that a {\it 
classical environment\/}  $E=({\cal P}_1,{\cal P}_2)$ is 
a pair of sets of measurable partitions of $\Omega$, and 
that ${\cal X}_i$ is the set of $S_i$-valued random 
variables that are measurable with respect to some 
partition in ${\cal P}_i$.   

In the discussion to follow we fix an environment 
$E=({\cal P}_1,{\cal P}_2)$

{\bf Definition 6.5.}  
Given \G and $E$ as above, we define a new game of 
private 
information $\G(E) $ as 
follows:
\itemitem{\bull}The information sets are the ${\cal 
A}_i$
\itemitem{\bull}The strategy sets are the ${\cal X}_i$
\itemitem{\bull}The payoff functions are given by  
$$P_i(A_1,A_2,X,Y)=\int_{S_1\times S_2}
P_i(A_1,A_2,x,y)d\mu_{X,Y}(x,y)$$
where $\mu_{X,Y}$ is the probability distribution on 
$S_1\times S_2$ induced by $(X,Y)$.

{\bf Remarks 6.6.}  Starting with the game of private 
information \G and the 
environment $E$, we can first  introduce contingent 
strategies, yielding the ordinary game  $\G^\#$ (see 
6.3), and then allow players to randomize, forming the 
game $\G^\#(E)$ (see 3.3). 

Alternatively, we can first allow players to randomize, 
forming the game of private information $\G(E)$ (see 
6.5) and then  introduce contingent strategies,  
forming the game $\G(E)^\#$ (see 6.3 again).  

These games are equivalent (in the sense of 3.1) to 
Kuhn's games of mixed and behavioral strategies ([K]) 
and are therefore equivalent to each other.  

In fact, more is true; they are isomorphic.
Explicitly: in $\G^\#(E)$, a 
strategy is an element of $\Hom(\Omega,\Hom({\cal 
A}_i,S_i))$, whereas in $\G(E)^\#$, a strategy is an 
element of $\Hom({\cal A}_i,\Hom(\Omega,S_i))$.  The 
strategy sets are carried back and forth to each other 
by the inverse bijections $$\Phi:\Hom(\Omega,\Hom({\cal 
A}_i,S_i))
\longleftrightarrow \Hom({\cal A}_i,\Hom(\Omega,S_i))
:\Psi$$
defined by
$((\Phi(f))(a))(\omega)=(f(\omega))(a)$ and 
$((\Psi(g))(\omega))(a)=(g(a))(\omega)$, and these 
bijections induce an isomorphism of games.

(It's straightforward to check that the maps $\Psi$ and 
$\Phi$ preserve the necessary measurability conditions.)

Thus we have:  

{\bf Theorem 6.7.} The games $\G^\#(E)$ and $\G(E)^\#$ 
are isomorphic. \fn{Readers of a certain mathematical 
bent will note that we've constructed not just an 
isomorphism but a {\it natural\/} isomorphism, which 
entitles us to think of these games as essentially ``the 
same''.  Other readers won't miss much if they ignore 
this footnote.}

{\bf Remarks 6.8.}  We will see in the next section that 
in the quantum context, this simple and natural 
isomorphism does not exist.  The proof of 6.7 no longer 
works, essentially because it's not possible to think of 
all the quantum observables as random variables 
originating in the same sample space.   Indeed, we will 
see in Theorems 7.4 and 7.6 that no analogue of Theorem 
6.7 can 
possibly hold in the quantum case.

{\bf Remark 6.9.}  Figure 1 below might be useful for 
keeping track of the notation in Remark 6.6 and Theorem 
6.7:

\vskip -2in

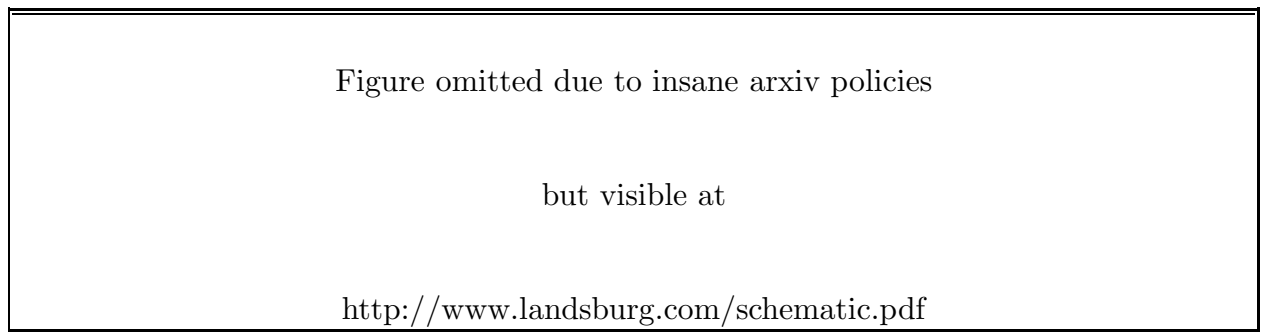
\vbox to 4in{
\centerline{\phantom{x}}
\boxit{
\centerline{Figure omitted due to insane arxiv policies}
\centerline{\phantom{x}}
\centerline{but visible at}
\centerline{\phantom{x}}
\centerline{http://www.landsburg.com/schematic.pdf}}}

(The section numbers indicate where these concepts are 
defined.)

\bigskip

\noindent{\bf 7.  Quantum Environments for Games of 
Private Information.}

Here we generalize as much of the previous section as 
possible to the quantum context.  We will see that in 
games of private information, there can be quantum 
equilibria that are in no sense equivalent to any 
classical correlated equilibrium.

We fix a game of private information \G.  For 
concreteness we take the strategy sets to be 
$S_1=S_2=\{\H,\T\}$.  We also fix a quantum environment 
$\xi$ as in 5.1.

{\bf Definition 7.1.}  (This is the quantum analogue of 
Definition 6.5).  Given \G and $\xi$ as above, we define 
a new game of private information $\G(\xi)$ as follows:
\itemitem{\bull}The signal sets are the ${\cal A}_i$
\itemitem{\bull}The strategy sets are the sets ${\cal 
U}_i$ defined in 5.2 and 5.3.
\itemitem{\bull}The payoff functions are given by 
$$P_i(A_1,A_2,U,V)=\int_{S_1\times 
S_2}P_i(A_1,A_2,s,t)d\mu_{UV}(s,t)$$
where $\mu_{UV}$ is the probability distribution defined 
in 5.4.

{\bf Remarks 7.2.}  Having formed the game of private 
information $\G(\xi)$, we can apply 6.3 to form the game 
$\G(\xi)^\#$, which is analogous to a game of behavioral 
strategies.  (In this game, a strategy consists of a map 
from the signal set $A_i$ to the set of quantum 
strategies.)  It then becomes natural to inquire about 
an analogue of Theorem 6.7.  

Therefore, let $(G,\xi)$ be a pair consisting of a game 
of private information, and a quantum environment.  
Figure Two is the quantum analogue of Figure One:

\vskip -2in

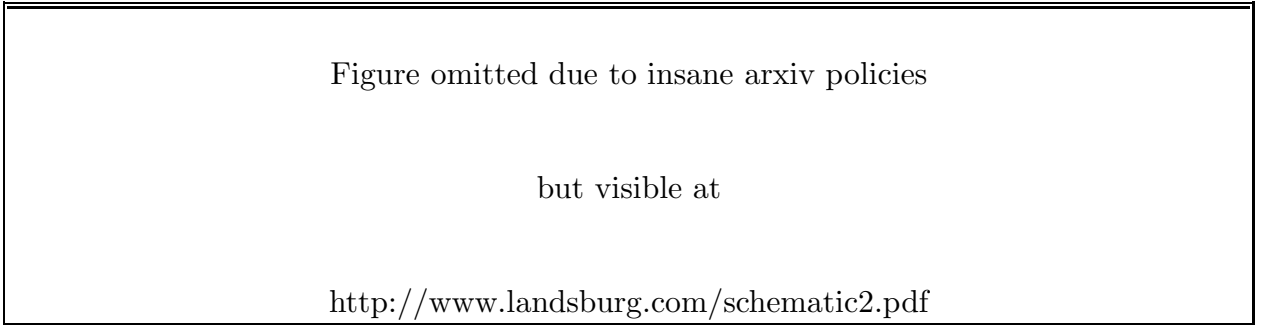
\vbox to 4in{
\centerline{\phantom{x}}
\boxit{
\centerline{Figure omitted due to insane arxiv policies}
\centerline{\phantom{x}}
\centerline{but visible at}
\centerline{\phantom{x}}
\centerline{http://www.landsburg.com/schematic2.pdf}}}

The first minor difficulty is that we haven't defined 
the bottom map.  That's because $\xi$, as defined in 
(5.1), is a quantum observable that takes only the two 
values \H and \T, and so serves as an environment only 
for games where the strategy sets $S_i$ have cardinality 
two.  In the game $\G^\#$, the strategy sets are 
typically much larger.

So if we want to define the bottom horizontal arrow, 
we'll need a new construction that converts the pair 
$(\G^\#,\xi)$ to a game $\G^\#(\xi)$.  It's not 
difficult 
to imagine reasonable constructions; we could, for 
example, allow players to observe multiple quantum coins 
so that the space of possible outcomes has the 
appropriate cardinality.

However, we don't have to worry about the details, 
because we're about to prove that {\it no\/} reasonable 
construction can work.

{\bf Definition 7.3.}  Let $\H$ be any (two-player but 
not necessarily two-by-two) game with 
with strategy sets $S_i$ and payoff functions $P_i$.  
Let $\H'$ be any other such game, with strategy sets 
$S_i'$ and payoff functions $P_i'$.  

We say that $\H'$ is a {\it stochastic extension} of 
\H if for any $(s',t')\in S_1'\times S_2'$, there is a 
probability distribution $\mu$ on $S_1\times S_2$ such 
that 
$$P_i'(s',t')=\int_{S_1\times S_2}P_i(s,t)d\mu(s,t)$$  

For example, if \G is any game, $E$ any classical 
environment, and $\xi$ any quantum environment, then 
$\G(E)$ and $\G(\xi)$ are both stochastic extensions of 
\G.

{\bf Theorem 7.4.}  There exists a game of private 
information $\G$ and a quantum environment $\xi$ such 
that $\G(\xi)^\#$ is not a stochastic extension of 
$\G^\#$.

{\bf Remark 7.4.1.}  Because the (still undefined) 
bottom arrow in Figure Two is supposed to be the 
analogue of ``allowing random (or quantum) strategies'', 
it is reasonable to require it to map any pair 
$(\H,\xi)$ to some stochastic extension of \H.  Theorem 
7.4 says that if we impose this requirement, then no 
matter how we define the bottom map, the analogue of 
Theorem 6.7 must fail.

{\bf Proof of 7.4.}  Let \G be the following game of 
private information:
\itemitem{a)}The strategy sets are $S_1=S_2=\{\H,\T\}$
\itemitem{b)}The signal sets are ${\cal A}_1={\cal A}_2=
\{\hbox{cat,dog}\}$
\itemitem{c)}The probability distribution on ${\cal 
A}_1\times {\cal A}_2$ is uniform
\itemitem{d)}The payoff functions are given by (2.1.1).

Let $\xi$ be the quantum environment $\HH+\TT$.  Then in 
the game $\G^\#$, the maximum possible payoff is .75.  
But it follows from example 5.9 that in the game 
$\G(\xi)^\#$, it is possible to achieve a payoff of 
approximately .85.  Therefore $\G(\xi)^\#$ cannot be a 
stochastic extension of $\G^\#(\xi)$.

{\bf Remark 7.5.}  Far more generally, if $\xi$ is any 
non-classical quantum environment, then there are (by 
definition) finite families ${\cal U}_1$ and ${\cal 
U}_2$ of unitary operators such that the various 
strategy pairs $(U_1,U_2)$ $(U_i\in{\cal U}_i)$ yield a 
family of probability distributions on $S_1\times S_2$ 
that cannot be mimicked by any classical random 
variables.  Thus if the signal sets in some game \G have 
cardinalities at least as large as those of the ${\cal 
U}_i$, then 
$\G(\xi)^\#$ cannot be  a stochastic extension of 
$\G^\#$.
Thus we can state:

{\bf Theorem 7.6.}  Let $\xi$ be a quantum environment, 
and suppose that for every game of private information 
\G, that $\G(\xi)^\#$ is a stochastic extension of 
$\G^\#(\xi)$.  Then $\xi$ is classical (as defined in 
5.8).  

Less formally:  Theorem 7.4 shows that there exists a 
quantum environment for which there is no reasonable
analogue of the Kuhnian construction for games
 of mixed strategies; 
Theorem 7.6 says that the same is 
true in any non-classical quantum environment.

{\bf Remark 7.7.}  If \G is a game of private 
information, then any quantum equilibrium in the 
game $\G^\#$  is a correlated equilibrium in $\G^\#$ and 
hence {\it a fortiori} (see 3.6) a Nash equilibrium in 
some stochastic extension of $\G^\#$.  

But if $\xi$ is a quantum environment, then $\G(\xi)^\#$ 
is {\it not\/} in general a stochastic extension of 
$\G^\#$, so it's plausible that it will have new 
equilibria which do not come from correlated equilibria 
in the game $\G^\#$.  In Section 9, we will see 
explicit examples of exactly this phenomenon.

\bigskip

\noindent{\bf 8.  Prelude to the Examples.}

In Section 9, we will offer two examples to illustrate 
the value of quantum strategies in games of private 
information.  To compute equilibria, we will exploit 
some special features of those examples, which we 
highlight here.

{\bf General Setup 8.1.}  In our examples, the strategy 
sets and signal sets will all have cardinality two.  We 
can take both signal sets to be $\{\C,\D\}$ and both 
strategy sets to be $\{\H,\T\}$.  Then a game of private 
information is specified by a probability distribution 
over the four signal pairs \CC, \CD, \DC, \DD, together 
with four two-by-two game matrices, one for each of 
these signal pairs.  In our examples, we will always 
take the last three of these four game matrices to be 
the same.  Thus the payoff structure is specified as 
follows:

$$\matrix
{\hbox{If both players receive signal \C}&&
\hbox{If either player receives signal \D}\cr
\matrix{
&\hbox{\bf Player Two}\cr
\matrix{\phantom{q}\cr
\matrix{ \vbox{\bf \hbox{Player}\hbox{One}}&
\matrix{}\cr}\cr}&
\matrix{
&{\H}&{\T}\cr
\H&(A,B)&(C,D)\cr
\T&(E,F)&(G,H) \cr}
\cr}
&&
\matrix{
&\hbox{\bf Player Two}\cr
\matrix{\phantom{q}\cr
\matrix{ \vbox{\bf \hbox{Player}\hbox{One}}&
\matrix{}\cr}\cr}&
\matrix{
&{\H}&{\T}\cr
\H&(I,J)&(K,L)\cr
\T&(M,N)&(P,Q) \cr}
\cr}\cr}\eqno(8.1.1)$$

Let \G be a game of private information with payoff 
structure (8.1.1) and let $\xi$ be the quantum 
environment $\HH+\TT$.  We will study the game of 
behavioral strategies $\G(\xi)^\#$.  Here a strategy for 
Player One is a pair $(U_{\C}, U_{\D})$ where each $U_i$ 
is either a special unitary matrix or one of the pure 
strategies \H and \T.   (A special unitary matrix 
corresponds to physically manipulating one's coin before 
observing it; a pure strategy involves throwing one's 
coin away.)

Likewise, a strategy for Player Two is a pair 
$(V_{\C},V_{\D})$.  

We will (temporarily) restrict attention to the subgame 
$\G(\xi)^\#_0$, in which $U_{\C},U_{\D},V_{\C}$ and 
$V_{\D}$ are restricted to be special unitary matrices
 --- that is, pure strategies are not allowed.   (This 
is unlikely to be a very interesting game in practice, 
but computing equilibria in $\G(\xi)^\#_0$ will be a 
stepping stone to computing equilibria in $\G(\xi)^\#$.) 

In the game $\G(\xi)^\#_0$, Example (4.6) implies that 
Player One's payoff is computed as follows: 

$$\eqalign
{& 
\Prob(\CC)\Bigg(s\Big(U_{\C}V^T_{\C}\Big){A+G\over2}+
t\Big(U_{\C}V^T_{\C}\Big)
{C+E\over2}\Bigg)\cr
+&
\Prob(\CD)\Bigg(s\Big(U_{\C}V^T_{\D}\Big){I+P\over2}+
t\Big(U_{\C}V^T_{\D}\Big)
{K+M\over2}\Bigg)\cr
+&
\Prob(\DC)\Bigg(s\Big(U_{\D}V^T_{\C}\Big){I+P\over2}+
t\Big(U_{\D}V^T_{\C}\Big)
{K+M\over2}\Bigg)\cr
+&
\Prob(\DD)\Bigg(s\Big(U_{\D}V^T_{\D}\Big){I+P\over2}+
t\Big(U_{\D}V^T_{\D}\Big)
{K+M\over2}\Bigg)\cr}\eqno(8.1.2)$$
where, for any special unitary matrix 
$U=\pmatrix{x&y\cr\-\overline{y}&\overline{x}\cr}$, we 
set $s(U)=|x|^2$ and $t(U)=|y|^2$.

We define the payoff structure (8.1.1) to be {\sl 
balanced\/} if the following four equations hold:

$$\matrix{
A+G=B+H&\qquad\qquad\qquad\qquad&C+E=D+F\cr
I+P=J+Q&&K+M=L+N\cr}\eqno(8.1.3)$$

If the original game \G is balanced, then in the quantum 
game $\G(\xi)^\#_0$ the two players' payoff functions 
are identical.  This implies:

{\bf Proposition 8.2.}  If \G is balanced, then there is 
a Pareto optimal equilbrium in the game $\G(\xi)^\#_0$.

{\bf Proof.}  Because the space of special unitary 
matrices is compact, there exist matrices 
$U_{\C},U_{\D},V_{\C},V_{\D}$ that maximize the value of 
(8.1.2).  

{\bf Discussion 8.3.}  If we assume that players achieve 
a Pareto optimal equilibrium in the game $\G(\xi)^\#_0$, 
then the discussion above reduces the computation of 
that equilibrium to a maximization problem over 4-tuples 
of special unitary matrices.  This can be 
computationally quite cumbersome.  
Our next result will demonstrate that we can restrict 
our attention to a small subset of those 4-tuples, which 
will make it possible to do these computations by hand.

{\bf Definition 8.4.}  For $\theta\in[0,2\pi)$, set (as 
in Example 5.9) 
$$M(\theta)=\pmatrix{\cos(\theta)&\sin(\theta)\cr
-\sin(
\theta)&\cos(\theta)}$$

{\bf Theorem 8.5.}  In the setup of (8.1), suppose that
\G is balanced.  Then there exists a real number
$\theta$ such that the maximum of (8.1.2) is achieved at
$$U_{\C}=M(-\theta)\quad U_{\D}=M(\theta)\quad
V_{\C}=M(2\theta) \quad V_{\D}=M(0)\eqno(8.5.1)$$

{\bf Proof.}  See the appendix to this paper.

{\bf Corollary 8.6.}  In the setup of (8.1), if $\G$ is
balanced, then in $\G(\xi)^\#_0$, there is a Pareto
optimal equilibrium of the form (8.5.1).

\bigskip

\noindent{\bf 9.  Examples.}

{\bf Example 9.1:  Cats and Dogs Revisited.}  In section 
2, we introduced the cat/dog game (adapted from [CHTW]) 
and claimed that players could improve their payoffs via 
quantum strategies.  In this section, we will return to 
that game and compute the equilibrium outcome when 
quantum strategies are available.  

As in Section 2, each player is asked ``Do you like
cats'' or ``Do you like dogs'', with the questions
chosen independently via fair coin flips.  Players then
answer yes or no, and receive the following payoffs:

$$\matrix
{\hbox{If both players are asked about cats}&&
\hbox{If either player is asked about dogs}\cr
\matrix{
&\hbox{\bf Player Two}\cr
\matrix{\phantom{q}\cr
\matrix{ \vbox{\bf \hbox{Player}\hbox{One}}&
\matrix{\cr}\cr}\cr}&
\matrix{
&\hbox{\bf YES}&\hbox{\bf NO}\cr
\hbox{\bf YES}&(0,0)&(1,1)\cr
\hbox{\bf NO}&(1,1)&(0,0) \cr}
\cr}
&&
\matrix{
&\hbox{\bf Player Two}\cr
\matrix{\phantom{q}\cr
\matrix{ \vbox{\bf \hbox{Player}\hbox{One}}&
\matrix{\cr}\cr}\cr}&
\matrix{
&\hbox{\bf YES}&\hbox{\bf NO}\cr
\hbox{\bf YES}&(1,1)&(0,0)\cr
\hbox{\bf NO}&(0,0)&(1,1) \cr}
\cr}\cr}$$

Note that this game has the form of (8.1.1) and is 
balanced
in the sense of (8.2).  Thus we can apply Corollary 8.6 
to conclude
that there is a Pareto-optimal equilibrium in 
$\G(\xi)^\#_0$ of the form (8.5.1):

$$U_{\C}=M(-\theta)\quad U_{\D}=M(\theta)\quad
V_{\C}=M(2\theta) \quad V_{\D}=M(0)$$

(Here \C and \D, of course, stand for ``cats'' and 
``dogs''.)

At this equilibrium, expression (8.1.2), for the common 
payoff to the two players, reduces to 
$${1\over4}\sin^2(3\theta)+{3\over4}\cos^2(\theta)\eqno(
9.1.1)$$
This expression is maximized at $\theta=\pi/8$, where it 
takes the approximate value .85.  This, then, is a 
Pareto optimal equilibrium in $\G(\xi)^\#_0$.

It is easy to verify that this remains a Pareto optimal 
equilibrium in the full quantum game $\G(\xi)^\#$ (that 
is, it remains both deviation-proof and Pareto-optimal 
when players are allowed to choose pure strategies).  In 
example 9.2, we will carry out such a verification in 
detail.     

{\bf Remark 9.1.1.}  The physical interpretation of this 
equilibrium is that Player One rotates his coin through 
an angle $-\pi/4$ or $\pi/4$ depending on whether he 
gets the cat question or the dog question, while Player 
Two rotates his coin through angle $\pi/2$ or $0$.
(Note that these are the same operations we considered 
in Example 5.9, where we used them to illustrate that no 
classical apparatus could mimic the resulting 
probability distribution.)

{\bf Remark 9.1.2.}  Not only do players win with 
probability .85 overall, they also win with probability 
.85 conditional on being asked any particular pair of 
questions.  

{\bf Remark 9.1.3.}  We have assumed here that players 
choose behavioral strategies --- maps from the signal 
set $\{\C,\D\}$ to the quantum strategy set $[0,2\pi)$.  
Suppose instead that players were required to choose 
(the quantum analogoue of) mixed strategies, so that 
each player throws a four-sided quantum coin yielding 
one of the four possible maps 
$\{\C,\D\}\rightarrow\{\hbox{Yes,No}\}$.  The result 
would be a correlated equilibrium in the four-by-four 
game with these strategies.  No such correlated 
equilibrium can yield an expected payoff greater than 
3/4.  This, reiterates the point  of Theorems 7.4 and 
7.5, i.e. the non-equivalence of 
mixed and behavioral strategies in the quantum context.

{\bf Remark 9.1.4.}  One might imagine the players 
pooling their information to achieve a better 
outcome.  It's important to recognize, however, that our 
quantum devices {\it in no way\/} allow such information 
pooling.  (Indeed, with information pooling, players 
could easily earn a certain payoff of 1.)  As we 
stressed in Remark 2.2, players do not communicate any 
aspect of their private information either to each other 
or to anyone else.

{\bf Example 9.2:  Airline Pricing.}   In Example 9.1, 
players have a 
shared goal.  Our next example illustrates similar 
phenomena in a game of greater economic interest, 
arising from a model of price competition with uncertain 
demand.  

Two identical airlines serve two types of customers.  
Low-demand customers have a reservation price $L$; high-
demand customers have a reservation price $H$.  There is 
a fixed population of $2x$ low-demand customers.  There 
is an uncertain population of high-demand customers.

First the airlines receive imperfect signals about the 
population of high-demand customers.  Then they set 
prices.  Then the high-demand customers, if any, arrive 
and buy seats.  Finally, the low-demand customers arrive 
and buy any remaining available seats at a low price.

The airlines' signals --- either $\N$ (negative) or $\P$ 
(postive) are drawn independently, with $\N$ and $\P$ 
equally probable.  If either signal is $\N$, there are 
$2y$ high-demand customers, for some $y<x$; otherwise 
there are none.

Each airline has a capacity constraint equal to the 
number of low-demand customers.  Their only costs are 
the fixed costs $F$ of running a flight.  This defines a 
game \G with the following payoff structure (where $\L$ 
and $\H$ stand for ``announce the low price $L$'' and 
``announce the high price $H$''):

$$\matrix
{\hbox{If both firms receive signal \N}&&
\hbox{If either firm receives signal \P}\cr
\matrix{
&\hbox{\bf Firm Two}\cr
\matrix{\phantom{q}\cr
\matrix{ \vbox{\bf \hbox{Firm}\hbox{One}}&
\matrix{\L\cr \H\cr}\cr}\cr}&
\matrix{
{\L}&{\H}\cr
(A,A)&(B,0)\cr
(0,B)&(0,0) \cr}
\cr}
&&
\matrix{
&\hbox{\bf Firm Two}\cr
\matrix{\phantom{q}\cr
\matrix{ \vbox{\bf \hbox{Firm}\hbox{One}}&
\matrix{\hbox{\L}\cr\hbox{\H}\cr}\cr}\cr}&
\matrix{
{\L}&{\H}\cr
(C,C)&(B,0)\cr
(0,B)&(D,D) \cr}
\cr}\cr}\eqno(9.2.1)$$

Here $A=xL-F$, $B=2xL-F$, $C=(x+y)L-F$, and $D=yH-F$.   
In
particular, $A<C<B$.  Note that if $D<B$ then \L is
always a dominant strategy for both players, so to keep
things interesting we will assume $B<D$.

The analysis of this game depends heavily on the values 
of
$A$, $B$, $C$ and $D$.  To avoid a proliferation of 
cases,
and to focus on a particularly interesting example, we 
take
$x=49$, $y=19$, $L=1$, $H=108/19\approx 5.68$, and 
$F=48$.
This gives payoffs of $A=1$, $B=50$, $C=20$, $D=60$.

\medskip

{\bf 9.2.2.  Classical Equilibrium.}
A (mixed) strategy for Firm One is a pair of 
probabilities
$(p_{\bf N},p_{\bf P})$, with $p_i$ the probability of 
playing $L$ after receiveing signal $i$.

Similarly,
Firm Two's strategy is a pair $(q_{\bf N},q_{\bf P})$.

When Firm One receives signal $\N$, it is equiprobable 
that Firm Two has received either signal $\N$ or signal 
$\P$.
Firm One's expected payoff is then
$$
{1\over2}\Big(p_{\bf N}q_{\bf N}\cdot 1+
              p_{\bf N}(1-q_{\bf N})\cdot 50\Big)+
{1\over2}\Big(p_{\bf N}q_{\bf P}\cdot 20+p_{\bf N}(1-
q_{\bf P})\cdot 50+
(1-p_{\bf N})(1-q_{\bf P})\cdot 60\Big)$$

When Firm One receives a positive signal, the expected
payoff is
$$
{1\over2}\Big(p_{\bf P}q_{\bf N}\cdot 20+
              p_{\bf P}(1-q_{\bf N})\cdot 50+(1-p_{\bf 
P})(1-q_{\bf N})60\Big)+
{1\over2}\Big(p_{\bf P}q_{\bf P}\cdot 20+p_{\bf P}(1-
q_{\bf P})\cdot 50+
(1-p_{\bf P})(1-q_{\bf P})\cdot 60\Big)$$

Firm One chooses $p_{\bf N}$ and $p_{\bf P}$ to maximize 
these
expressions and Firm Two behaves symmetrically.
One checks that the unique equilibrium is at
$p_{\bf N}=p_{\bf P}=q_{\bf N}=q_{\bf P}=1$; that is, 
both firms always play
$\L$.  This guarantees them payoffs of 1 when both
receive
signal $\N$ (1/4 of the time) and 20 when either 
receives
a
signal $\P$ (3/4 of the time), for an expected payoff of
$15.25$.

{\bf 9.2.3.  Quantum Equilibrium.}  Next we equip our 
players with a pair of coins in the state $\xi=\HH+\TT$.  
Because the payoff structure 
(9.2.1) satisfies the assumptions  of Theorem (8.5), the 
game $\G(\xi)^\#_0$ (that is, the game in which players 
are required to adopt quantum strategies) has a Pareto 
optimal equilibrium of the form (8.5.1).   

At this equilibrium, when both players receive signal 
\N, their payoff is 
$$\cos^2(3\theta){A\over 
2}+\sin^2(3\theta){B\over2}\eqno(9.2.3.1)$$
and otherwise their payoff is
$$\cos^2(\theta){C+D\over2}+\sin^2(\theta){B\over2}\eqno
(9.2.3.2)$$
Adding (9.2.3.1) to 3 times (9.2.3.2) and substituting 
the assumed values for $A,B,C$, and $D$, the payoff 
becomes

$$\cos^2(3\theta){1\over2}+\sin^2(3\theta){50\over2}
+3\cos^2(\theta){20+60\over2}+3\sin^2(\theta)
{50\over2}$$

It is easy to verify that this expression is maximized 
at 

$$\theta=\hbox{ArcCos}\left({1\over2}\sqrt{
{14+\sqrt{79}\over7}}\right)\eqno(9.2.3.3)$$
where it takes the value
$${3087+79\sqrt{79}\over112}\approx 33.83$$
This beats the classical payoff of 15.25.

Thus we have found a quantum equilibrium in 
$\G(\xi)^\#_0$.  We want to show that it remains an 
equilibrium in the full quantum game $\G(\xi)^\#$, where 
players are allowed to adopt pure strategies.  That is, 
we want to check that neither player wants to 
deviate from the quantum equilibrium by playing 
classically.

Write $({\bf s}_1, {\bf s}_2)$ for the strategy pair 
described by (8.5.1) and (9.2.3.3)

{\bf  Claim 1:} If Firm Two plays the quantum strategy 
$\s_2$,
and if Firm One receives the negative signal $N$, then 
Firm
One
plays the quantum strategy $\s_1$.

{\bf Proof.}  Note that Firm Two plays $\L$ and $\H$ 
with equal
probability.  Therefore, when Firm One receives a 
negative
signal, it can earn any of the following 
returns, with equal probability:

$$\matrix{
\hbox{Strategy}&\hbox{Return}\cr
\phantom{duh}\cr
\L&{1\over4}(1+50+20+50)=30.25\cr
\H&{1\over4}(0+0+0+60)=15\cr
{\bf s_1}&{181+7\sqrt{79}\over8}\approx30.40\cr}$$

Because $30.40>30.25$ (and because the quantum
strategy ${\bf s}$ maximizes both players' payoffs over 
all
alternative quantum strategies), Firm One chooses 
strategy
${\s_1}$.

{\bf  Claim 2:} If Firm Two plays the quantum strategy 
$\s_2$,
and if Firm One receives the positive signal $P$, then 
Firm
One
plays the quantum strategy $\s_1$.

{\bf Proof.}  In this case, Firm One's payoffs are
known to be 

$$\matrix{
\hbox{Strategy}&\hbox{Return}\cr
\phantom{duh}\cr
\L&{1\over2}(20+50)=35\cr
\H&{1\over2}(0+60)=30\cr
{\bf 
s_1}&{65\over2}+{15\sqrt{79}\over28}\approx37.26\cr}$$
The result follows 
because
$37.26>35$.

Combining Claims 1 and 2, we see that when Firm Two 
plays $\s_2$, Firm One does not want to deviate; in view 
of the symmetries of the game, the same is of course 
true with the firms reversed.   
Thus $(\s_1,\s_2)$ is genuinely an equilibrium, 
and,
as
we have already seen, it is Pareto superior to the 
unique
classical equilibrium where both firms always play $\L$.

\medskip

{\bf 9.2.3.  Correlated Equilibrium.}  
Consider the game in which a strategy is a behavioral 
strategy, i.e. a map from the 
set 
of signals $\{\N,\P\}$ to the set of actions 
$\{\L,\H\}$.  
The 
conclusion 
of section (9.2) is equivalent to the statement that 
the 
only 
mixed strategy equilibrium in this game is the pair 
$(1_{\L},1_{\L})$ where $1_{\L}$ is the constant map 
``always 
play ${\L}$''.  It's not hard to check that this is also 
the 
only correlated equilibrium in the sense of Aumann [A].  

In other words, in a world governed by the usual laws of 
probability theory, players who condition their 
strategies on 
both the signals they receive and their observations of 
(possibly correlated) random variables can still do no 
better than in the classical equilibrium of section 
(9.2.1).  
In still other words, a referee who dictates conditional 
strategies (subject to a deviation-proofness criterion) 
cannot improve the outcome.  
This makes it all the more striking that they {\it 
can\/} 
improve the outcome in the quantum world of section 
(9.2.2).

{\bf 9.2.4.  Welfare.}  Consumer surplus occurs only 
when 
high
demand customers pay low prices.  In any of these cases, 
38 high-demand
customers earn surpluses of 89/19, for a total surplus 
of
178.

In classical equilibrium, both firms play $\L$, ending 
up 
in
the upper left corner of the left-hand payoff matrix 1/4 
of
the time and the upper left corner of the right-hand 
payoff
matrix 3/4 of the time.  Thus producer surplus is
$(1/4)\cdot 2+(3/4)
\cdot 40=30.5$ and consumer surplus 
is
$(3/4)\cdot 178=133.5$.

If firms could collude, they would play either $(\L,\H)$ 
or
$(\H,\L)$ in the low-demand state of the world and 
$(\H,\H$)
in the high demand state of the world for a producer 
surplus
of $(1/4)\cdot 50+(3/4)
\cdot 120=102.5$, while consumer
surplus would fall to zero.  Thus it makes sense that a
regulator would want to prohibit collusion.

In quantum equilibrium, one computes that producer 
surplus
is approximately 67.66 and consumer surplus is 
approximately
78.94.

In summary:
$$\matrix{&\hbox{Classical}&\quad\hbox{Quantum}&\quad
\hbox{Classical with Collusion}\cr
\hbox{Consumer Surplus}&133.5&78.94&0\cr
\hbox{Producer Surplus}&30.5&67.66&102.5
\cr
\hbox{Total}&164&146.6&102.5\cr}$$

Thus our regulator would want to prohibit the use of 
quantum
technology, though not as much as he wants to prohibit
collusion (and firms would want to use quantum 
technology,
though not as much as they want to collude).  The 
quantum
technology, however, is quite undetectable (after a firm
announces it's strategy, how do you know whether it
randomized by looking at a classical or a quantum 
coin?), and hence presumably impossible to regulate.

\bigskip

\centerline{\bf Appendix:  Proof of Theorem 8.5.}

Given a special unitary matrix 
$$U=\pmatrix{P&Q\cr-\overline{Q}&\overline{P}}$$
set $t(U)=|Q|^2$. 

Given four special unitary matrices $A,B,S,T$, set
$$X=t(AS)\quad Y=t(BS) \quad Z=t(BT)\quad 
W=t(AT)\eqno(A.1)$$

Let $p,q$ and $r$ be any real numbers.

Consider the following function

$$f(A,B,S,T)=pW+q(X+Y+Z)+r\eqno(A.2)$$

{\bf Theorem A.3.} There exists a real number $\theta$ 
such that the maximum of (A.2) is achieved at
$$A=M(\theta)\quad B=M(-\theta)\quad S=M(0) \quad 
T=M(2\theta)\eqno(A.4)$$

At this maximum we have 

$$X=Y=Z=\sin^2(\theta)\qquad 
W=\sin^2(3\theta)\eqno(A.5)$$

{\bf Proof.}  From [Land] (unnumbered proposition on 
page 455) it follows that $(X,Y,Z,W)\in[0,1]^4$ is in 
the image of $f$ if and only if 
$$|XY-X-Y-ZW+Z+W| \le 2(\sqrt{X-X^2}\sqrt{Y-
Y^2}+\sqrt{Z-Z^2}\sqrt{W-W^2})\eqno(A.6)$$

We seek to maximize $p(X+Y+Z)+qW$ over $[0,1]^4$ subject 
to the constraint (A.6).  Note that if $(X,Y,Z,W)$ 
satisfies the constraint, then so does $(M,M,M,W)$ where 
$M$ is the average of $X,Y$ and $Z$.  Thus to find a 
maximum, we can assume that $X=Y=Z$.  This reduces the 
problem to a constrained maximization over two variables 
$X$ and $W$ which can be solved (somewhat laboriously) 
by hand.  This gives $X$ and $W$ as explicit functions 
of $p$ and $q$; by inspection,  we find that 
$\hbox{ArcSin}(\sqrt{W})=3\hbox{ArcSin}(\sqrt{X})$; thus 
we can 
set
$X=\sin^2(\theta)$ and $W=\sin^2(3\theta)$.  Therefore 
(A.5) 
is a sufficient condition for a maximum, and this 
condition holds at (A.4). 

{\bf Corollary A.7:}  Theorem 8.5 is true.

{\bf Proof.}  Note that expression (8.1.2) is of the 
form A.2 (use the fact that s(x) = 1-t(x)).  Therefore 
Theorem A.3 applies.

\bigskip

\centerline{\bf References}

\item{[A]}Aumann, ``Subjectivity and Correlation in
Randomized Strategies'', {\it J. Math Econ\/} 1 (1974).

\item{[B]}Adam Brandenburger, ``The relationship between
quantum and classical correlation in games'', {\it Games
and Economic Behavior} 69 (2010), 175-183.

\item{[CHTW]}R.~Cleve, P.~Hoyer, B.~Toner, and
J.~Watrous,
``Consequences and Limits of Nonlocal Strategies'', {\it
Proc. of the 19th Annual Conference on Computational
Complexity\/} (2004), 236-249

\item{[EW]}J.~Eisert and M.~Wilkens, ``Quantum Games'',
{\it J.
of Modern Optics\/} 47 (2000), 2543-2556

\item{[EWL]} J. Eisert, M. Wilkens and
M. Lewenstein, ``Quantum Games and
Quantum Strategies'', {\it Phys. Rev.
Letters\/} 83, 3077 (1999).

\item{[K]}H. Kuhn, ``Extensive Games and the Problem of
Information'', in {\it Contributions to the Theory of
Games
II}, Annals of Math Studies 28 (1953).

\item{[L]}S.~Landsburg, ``Nash Equilibria in Quantum
Games'', to appear in {\it Proceedings of the American
Mathematical Society\/}, December 2011.

\item{[LaM]}
La Mura, P., ``Correlated Equilibria of Classical Strategic Games with
Quantum Signals,'' {\it International Journal of Quantum Information\/}, 3,
2005, 183-188.

\item{[Land]}L.J.~Landau, ``Empirical Two-Point 
Correlation Functions'', {\it Foundations of Physics} 18 
(1988).

\item{[Le]}D. Levine, ``Quantum Games Have No News for
Economists'', working paper.

\bye